\documentclass[letterpaper, 10 pt, conference]{ieeeconf}  

\IEEEoverridecommandlockouts

\usepackage{cite}
\usepackage{amsmath,amssymb,amsfonts}
\usepackage{graphicx}
\usepackage{textcomp}
\def\BibTeX{{\rm B\kern-.05em{\sc i\kern-.025em b}\kern-.08em
    T\kern-.1667em\lower.7ex\hbox{E}\kern-.125emX}}

\usepackage{xcolor}
\usepackage{bm}

\usepackage{algpseudocode}
\usepackage{algorithm}

\newcounter{defcounter}
\newenvironment{pequation}{%
\addtocounter{equation}{-1}
\refstepcounter{defcounter}
\renewcommand\theequation{{\bf{P\thedefcounter}}}
\begin{equation}}
{\end{equation}}

\usepackage{subfigure}

\usepackage{cases}
\usepackage{array}

\usepackage{empheq}

\usepackage{dsfont}

\newcommand{\mybold}{\text{\usefont{U}{bbold}{m}{n}1}}
\MakeRobust{\mybold}

\usepackage{subfig}

\begin{document}

\title{\LARGE \bf On Privacy Preservation of Electric Vehicle Charging Control via State Obfuscation}

\author{Xiang Huo$^{\dagger}$ and Mingxi Liu$^{\dagger}$
\thanks{$^{\dagger}$X. Huo and M. Liu are with the Department of Electrical and Computer Engineering at the University of Utah, 50 S Central Campus Drive, Salt Lake City, UT, 84112, USA {\tt \small{\{xiang.huo, mingxi.liu\}@utah.edu}}.}
 \thanks{This work is supported by NSF Award: ECCS-2145408.}}


\maketitle
\thispagestyle{empty}
\pagestyle{empty}

\begin{abstract}

The electric vehicle (EV) industry is rapidly evolving owing to advancements in smart grid technologies and charging control strategies. While EVs are promising in decarbonizing the transportation system and providing grid services, their widespread adoption has led to notable and erratic load injections that can disrupt the normal operation of power grid.  Additionally, the unprotected collection and utilization of personal information during the EV charging process cause prevalent privacy issues. To address the scalability and data confidentiality in large-scale EV charging control, we propose a novel decentralized privacy-preserving EV charging control algorithm via state obfuscation that 1) is scalable \emph{w.r.t.} the number of EVs and ensures optimal EV charging solutions; 2) achieves privacy preservation in the presence of honest-but-curious adversaries and eavesdroppers; and 3) is applicable to eliminate privacy concerns for general multi-agent optimization problems in large-scale cyber-physical systems. The EV charging control is structured as a constrained optimization problem with coupled objectives and constraints, then solved in a decentralized fashion. Privacy analyses and simulations demonstrate the efficiency and efficacy of the proposed approach.


\end{abstract}


\section{Introduction}

The ongoing advancements in electric vehicle (EV) technologies have accelerated the development of a sustainable power grid, owing to the EVs' green credentials and flexible charging options. Despite the multifarious benefits, the occurrence of plug-and-play EV charging events, especially those involving a significant number of EVs, can cause several negative impacts on the power grid, such as load profile fluctuations, voltage deviations, and increased power loss \cite{yilmaz2012review}.  Therefore, advancing scalable EV charging coordination and control strategies is of paramount importance to alleviate the strain on the power grid and ultimately provide synergistic grid-edge services, such as valley-filling, peak-shaving, and frequency regulation.

In enabling such synergy between EVs and the power grid, the EV charging control problem can be framed as a constrained optimization problem. Let $\bm{x}_{\hat{\imath}} = {[x_{\hat{\imath}}(1), \ldots, x_{\hat{\imath}}(T)]}^{\mathsf{T}}$ denote the charging profile of EV $\hat{\imath}$ during $T$ consecutive time slots, $\mathcal{X}_{\hat{\imath}}$ denote the local feasible set that contains the charging requirements of EV $\hat{\imath}$, and function $g(\cdot)$ denote the networked constraint function. Then, the EV charging control problem can be formulated into a constrained optimization problem as 
\begin{pequation} \label{sss1}
\begin{aligned}
&\text{min} \quad J(\left\{\bm{x}_{\hat{\imath}}\right\}_{\hat{\imath}=1}^{\hat{n}}) = F \left(\bm{p}_{b}, \bm{x}_{1}\ldots, \bm{x}_{\hat{n}}\right) \\
& \,\, \text{s.t.}  \, \quad  \bm{x}_{\hat{\imath}} \in \mathcal{X}_{\hat{\imath}}, ~ \forall \hat{\imath}=1,2, \ldots, \hat{n} \\
&\,\,\,\, \qquad 
g(\bm{x})  \leq  0
\end{aligned}
\end{pequation}where the cost function $F(\cdot) : \mathbb{R}^T \mapsto \mathbb{R}$ is assumed convex and differentiable,  $\bm{p}_{b} \in \mathbb{R}^T$ captures the baseline load of the network, $\bm{x} = {[\bm{x}_{\hat{\imath}}^{\mathsf{T}}, \ldots, \bm{x}_{\hat{n}}^{\mathsf{T}}]}^{\mathsf{T}}$, and $\hat{n}$ denotes the total number of EVs.

The use of scalable optimization methods, such as distributed and decentralized approaches, has gained popularity in solving \eqref{sss1}. In \cite{karfopoulos2012multi}, a distributed multi-agent EV charging control method was developed based on the Nash certainty equivalence principle to account for network impacts. Gan \emph{et al.} in \cite{gan2012optimal}
proposed a decentralized EV charging control algorithm with the objective of addressing the valley-filling problem using EVs' charging loads. To scale with the EV fleet size and the length of control periods, decentralized EV charging protocols were developed in \cite{zhang2016scalable} for network-constrained EV charging problems. In \cite{liu2017decentralized}, a decentralized EV charging control scheme was developed to achieve valley filling, meanwhile accommodating individual charging needs and distribution network constraints. To further improve the scalability,  a distributed optimization framework was proposed in \cite{huo2022two} to offer two-facet scalability over both the agent population size and network dimension. 


Besides scalability, the increased risk of privacy exposure is another major obstacle in deploying large-scale EV charging control strategies. To address the pressing need for privacy preservation in both EV charging control and generic multi-agent systems, one potential solution is using differential privacy (DP). 
Fiore and Russo in \cite{fiore2019resilient} designed a DP-based consensus algorithm for multi-agent systems where a subset of agents are adversaries. In \cite{nozari2016differentially}, a distributed functional perturbation framework was developed based on DP  to protect each agent's private objective function. In \cite{wang2022differentially}, DP-based distributed algorithms were designed to preserve privacy in finding the Nash equilibrium of stochastic aggregative games. Although DP-based methods are commonly adopted for privacy preservation, the inevitable trade-off between accuracy and privacy remains a major challenge in practical implementation.

Another frequently utilized method for preserving privacy involves cryptographic techniques, such as the Paillier cryptosystem and Shamir's secret sharing (SSS). In \cite{fang2021secure}, a Paillier-based privacy-preserving algorithm was proposed for securing the average consensus of networked systems with high-order dynamics. Zhang \emph{et al.} in \cite{zhang2021privacy} developed a privacy-preserving power exchange service system that uses data encryption to protect EV users' privacy.  In \cite{huo2021encrypted}, a decentralized privacy-preserving multi-agent cooperative optimization paradigm was designed based on cryptography for large-scale industrial cyber-physical systems. In \cite{zhang2018enabling}, a novel decentralized  privacy preservation approach was designed by integrating a partially homomorphic cryptosystem into the decentralized optimization architecture. Compared to encryption-based methods that rely on large integer calculations, SSS-based privacy-preserving approaches are more efficient in the computation of shares while offering information-theoretical security \cite{huo2022secret}. In \cite{li2019privacy}, an SSS-based privacy-preserving algorithm was developed to solve the consensus problem while concurrently protecting each individual's private information. Rottondi \emph{et al.} in \cite{rottondi2014enabling} designed a privacy-preserving vehicle-to-grid architecture based on SSS to ensure the confidentiality of the private information of EV owners from aggregators. Huo and Liu in \cite{huo2022distributed} proposed an SSS-based privacy-preserving EV charging control protocol, which eliminates the need for a system operator (SO) in achieving overnight valley filling. While cryptographic methods effectively achieve high levels of accuracy and privacy, the accompanying increased computation and communication complexity become the bottleneck in their practical use.  Non-cryptographic approaches like state decomposition (SD) decompose the true state into two sub-states, and only one sub-state is visible to others, therefore protecting the true value of the original state. However, state-of-the-art SD-based strategies are not applicable to solve \eqref{sss1} as they mainly focus on consensus problems \cite{wang2019privacy,zhang2022privacy}.

This paper aims to design a decentralized privacy-preserving optimization algorithm, which is  scalable and low in complexity, suitable for large-scale multi-agent optimization, specifically for EV charging control. The contributions of this paper are three-fold: 1) the proposed decentralized privacy-preserving algorithm can scale with the number of EVs and provide optimal decentralized EV charging solutions; 2) privacy preservation is achieved in the presence of honest-but-curious adversaries and external eavesdroppers; and 3) the proposed approach has low computing and communication overhead, making it widely applicable for preserving privacy in coupled multi-agent optimization problems in cyber-physical systems.



\section{Problem Formulation}

\subsection{Distribution Network Model}\label{Distribution-network}

In a radial distribution network, the power flow can be represented by DistFlow branch equations that consist of the real power, reactive power, and voltage magnitude \cite{baran1989network}. Consider a re-indexed radial distribution network and define $\mathbb{N}=\{i \mid i=1, \ldots, n\}$ as the set of downstream buses. Let $\left|V_{i}(t)\right|$ denote the voltage magnitude of bus $i$ at time $t$, $|V_0|$ denote the voltage magnitude of the slack bus, and $p_i(t) $ and $q_i(t)$ denote the active and reactive loads of bus $i$ at time $t$. 

\begin{figure}[!htb]
    \centering
\includegraphics[width=0.47\textwidth, trim={0.1cm 0cm 0.0cm 0cm},clip]{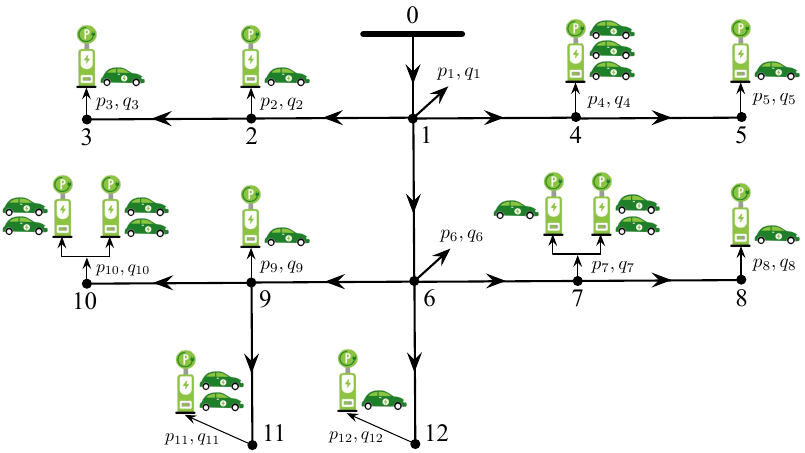}
    \caption{A 13-bus distribution network connected with EVs.}
\label{13_distribution_network}
     \vspace{-6mm}
\end{figure}

Following the linear DistFlow branch equations  \cite{baran1989network}, the squared voltage magnitude at node $i$ is
\begin{equation} \label{11sss}
\boldsymbol{V}_i =\boldsymbol{V}_{0}- 2\sum_{j=1}^{n} \boldsymbol{R}_{ij} \boldsymbol{p}_j - 2\sum_{j=1}^{n} \boldsymbol{X}_{ij} \boldsymbol{q}_j
\end{equation}
where $
\boldsymbol{V}_{i} =[\left|V_{i}(1)\right|^{2}, \ldots, \left|V_{i}(T)\right|^{2}]^{\mathsf{T}} \in \mathbb{R}^{T}$, $
\boldsymbol{V}_{0} =[\left|V_{0}\right|^{2}, \ldots, \left|V_{0}\right|^{2}]^{\mathsf{T}} \in \mathbb{R}^{T}$, $\bm{p}_i = [p_{i}(1), \ldots, p_{i}(T)]^{\mathsf{T}} \in \mathbb{R}^{T} $, $
\bm{q}_i = \left[q_{i}(1),\ldots, q_{i}(T)\right]^{\mathsf{T}} \in \mathbb{R}^{T}$, 
and the adjacency matrices $\bm{R}$ and $\bm{X}$ are defined as
\begin{align}
\boldsymbol{R} &\in \mathbb{R}^{n \times n},~  \boldsymbol{R}_{ij}=\sum_{(\hat{\imath}, \hat{\jmath}) \in \mathbb{E}_{i} \cap \mathbb{E}_{j}} r_{\hat{\imath} \hat{\jmath}} \nonumber\\ 
\boldsymbol{X} &\in \mathbb{R}^{n \times n}, ~\boldsymbol{X}_{ij}=\sum_{(\hat{\imath}, \hat{\jmath}) \in \mathbb{E}_{i} \cap \mathbb{E}_{j}} x_{\hat{\imath} \hat{\jmath}}\nonumber
\end{align} 
where $r_{\hat{\imath} \hat{\jmath}}$ and $x_{\hat{\imath} \hat{\jmath}}$ denote the  resistance and reactance from bus $\hat{\imath}$ to bus $\hat{\jmath}$, respectively. The sets of line segments that connect the slack bus to bus $i$ and bus $j$ are denoted by $\mathbb{E}_{i}$ and $\mathbb{E}_{j}$, respectively. In this paper, we focus on the charging control of plug-in EVs on radial distribution networks. A 13-bus distribution network with charging stations situated at different nodes is shown in Fig. \ref{13_distribution_network}.

\subsection{EV Charging Model}\label{EV-model}

Let $\bm{r}_{\hat{\imath}} \in \mathbb{R}^{T}$ denote the piece-wise constant charging power of the $\hat{\imath}$th EV during $T$ time intervals, and $\bm{r}_{\hat{\imath}}$ is constrained by 
\begin{equation}\label{1s}
    \bm{0} \leq \bm{r}_{\hat{\imath}} \leq \overline{\bm{r}}_{\hat{\imath}}
\end{equation}
where 
$\overline{\bm{r}}_{\hat{\imath}} = [\overline{r}_{\hat{\imath}}, \ldots, \overline{r}_{\hat{\imath}}]^{\mathsf{T}} \in \mathbb{R}^{T}$ and $\overline{r}_{\hat{\imath}}$ denotes the maximum charging power. 

Let $\delta_t$ denote the sampling period and $[1:T\delta_t]$ denote the charging duration.  To ensure EVs are charged to the desired energy levels by the end of the charging period, the total energy charged for  the $\hat{\imath}$th EV should satisfy
\begin{equation}
\bm{G}\bm{r}_{\hat{\imath}} = d_{\hat{\imath}}\label{2s}
\end{equation}
where $\bm{G} {=} [\delta_t\eta,\ldots,\delta_t\eta] {\in} \mathbb{R}^{1\times T}$ denotes the aggregation vector, $\eta$ denotes the charging efficiency, and $d_{\hat{\imath}}$ denotes the charging demand request of the $\hat{\imath}$th EV.

\subsection{Valley-Filling Optimization Problem}\label{opt-problem}

The valley-filling problem aims at filling the aggregated load valley and smoothing the aggregated load profile of the entire distribution network. This service is typically provisioned during the late evening and early morning when significant energy use reduction occurs. In this paper, the controllable charging loads of EVs, e.g., from community overnight parking and charging lot, are scheduled and shifted to flatten the valley in the power load profile. To this end, the valley-filling problem is formulated as a constrained optimization problem at the grid scale aiming at determining the optimal charging schedules of all EVs. 

Suppose in total $\hat{n}$ EVs need to be fully charged during the time period $[1:T\delta_t]$. This paper takes the nodal voltage constraint, which manifests as the global constraint, for example, to illustrate the impacts of EV charging on the distribution network as 
\begin{equation} \label{4sk}
    \underline{\bm{V}}  \leq  \bm{V}_i \leq \overline{\bm{V}}, ~ \forall i=1,2, \ldots, n 
\end{equation}
where $\underline{\bm{V}} {=} [\underline{V}, \ldots, \underline{V}]^{\mathsf{T}} {\in} \mathbb{R}^{T}$, $\overline{\bm{V}} {=} [\overline{V}, \ldots, \overline{V}]^{\mathsf{T}} {\in} \mathbb{R}^{T}$, $\underline{V}$ and $\overline{V}$ denote the lower and upper voltage bounds, respectively. 

The optimal EV charging control problem is then formulated into a quadratic programming problem as 
\begin{pequation} \label{s1}
\begin{aligned}
    &\underset{}{\text{min}} \quad J(\left\{\bm{r}_{\hat{\imath}}\right\}_{\hat{\imath}=1}^{\hat{n}}) = \frac{1}{2}\left\|\bm{p}_{b}+\sum_{\hat{\imath}=1}^{\hat{n}} \bm{r}_{\hat{\imath}} \right\|_{2}^{2} \\
& \,\, \text{s.t.}  \, \quad  \bm{r}_{\hat{\imath}} \in \mathcal{R}_{\hat{\imath}}, ~ \forall \hat{\imath}=1,2, \ldots, \hat{n} \\
& \,\,\,\, \qquad 
\underline{\bm{V}}  \leq  \bm{V}_i \leq \overline{\bm{V}}, ~ \forall i=1,2, \ldots, n
\end{aligned}
\end{pequation}where $\bm{p}_{b}$ denotes the aggregated baseline load and $\mathcal{R}_{\hat{\imath}}$ denotes the local feasible set of EV $\hat{\imath}$ that is defined by 
\begin{equation}
\mathcal{R}_{\hat{\imath}} = \{ \bm{r}_{\hat{\imath}} | \: 0 \leq \bm{r}_{\hat{\imath}} \leq \overline{\bm{r}}_{\hat{\imath}},   \bm{G}\bm{r}_{\hat{\imath}} = d_{\hat{\imath}} \}.
\end{equation}
Note that \eqref{1s} and \eqref{2s} are basic constraints that describe the EV charging process, additional constraints that introduce EVs' local characteristics can be included in the feasible set $\mathcal{R}_{\hat{\imath}}$ without affecting the algorithm design.

\section{Main Results}

\subsection{Decentralized PGM}

To solve the constrained optimization problem in \eqref{s1} via a decentralized manner, EVs (agents) can work cooperatively by adopting the projected gradient method (PGM) \cite{bertsekas2015parallel}. In PGM, the  $\hat{\imath}$th EV can update its decision variable (primal variable) by 
\begin{equation} \label{eq1}
\bm{r}_{\hat{\imath}}^{(\ell+1)} = \Pi_{\mathcal{R}_{\hat{\imath}}}[\bm{r}_{\hat{\imath}}^{(\ell)} - \gamma^{(\ell)}_{\hat{\imath}} \Phi_{\hat{\imath}}^{(\ell)} (\bm{r}^{(\ell)})] 
\end{equation}
where $\ell$ denotes the iteration index, $\bm{r}^{(\ell)} = [{\bm{r}_1^{(\ell)}}^\mathsf{T},\ldots,{\bm{r}_{\hat{n}}^{(\ell)}}^\mathsf{T}]^{\mathsf{T}}$, $\gamma^{(\ell)}_{\hat{\imath}}$ denotes the primal update step size of EV $\hat{\imath}$, $\Phi_{\hat{\imath}}^{(\ell)}(\cdot)$ denotes the first-order gradient of the Lagrangian function \emph{w.r.t.} $  \bm{r}_{\hat{\imath}}^{(\ell)}$, and $\Pi_{\mathcal{R}_{\hat{\imath}}}[\cdot]$ denotes the Euclidean projection operation onto $\mathcal{R}_{\hat{\imath}}$.

The relaxed Lagrangian of  \eqref{s1} can be derived as 
\begin{equation} \label{5}
    \mathcal{L}(\bm{r},\bm{\lambda})= \frac{1}{2}\left\|\bm{p}_{b}+\sum_{\hat{\imath}=1}^{\hat{n}} \bm{r}_{\hat{\imath}}\right\|_{2}^{2} + \sum_{i=1}^{n}\bm{\lambda}_i^\mathsf{T}(\underline{\bm{V}} - \bm{V}_i)
\end{equation}
where  $\bm{\lambda} = [\bm{\lambda}_1^\mathsf{T},\ldots,\bm{\lambda}_n^\mathsf{T}]^\mathsf{T}$ and $\bm{\lambda}_i$ denotes the dual variable associated with the $i$th inequality constraint. Note that the Lagrangian in \eqref{5} is relaxed by moving EVs' local constraints into $\mathcal{R}_{\hat{\imath}}$. Only the lower bound constraint on the bus voltage magnitudes is considered, as the charging loads of EVs are the only active power consumption within the distribution network.

The subgradients of $\mathcal{L}(\bm{r},\bm{\lambda})$ \emph{w.r.t.} $\bm{r}_{\hat{\imath}}$ and $\bm{\lambda}_i$ are
\begin{subequations}\label{6s}
\begin{align}
    \nabla_{\bm{r}_{\hat{\imath}} } \mathcal{L}(\bm{r},\bm{\lambda}) & =  \bm{p}_{b} + \sum_{\hat{\imath}=1}^{\hat{n}} \bm{r}_{\hat{\imath}}  -\sum_{i=1}^{n}\nabla_{\bm{r}_{\hat{\imath}} }(\bm{\lambda}_i^\mathsf{T}\bm{V}_i) \label{4as}\\
     \nabla_{\bm{\lambda}_i} \mathcal{L}(\bm{r},\bm{\lambda}) & = \underline{\bm{V}} - \bm{V}_i.   \label{4bs}
\end{align}
\end{subequations}
Substitute the linear DistFlow branch equation \eqref{11sss} into \eqref{6s}, we have 
\begin{subequations}\label{7s}
\begin{align}
    \nabla_{\bm{r}_{\hat{\imath}} } \mathcal{L}(\bm{r},\bm{\lambda}) & =  \bm{p}_{b} + \sum_{i=1}^{n} \bm{p}_{i} -  \hat{\bm{s}}_{\hat{\imath}} \label{7as}\\
     \nabla_{\bm{\lambda}_i} \mathcal{L}(\bm{r},\bm{\lambda}) & = \tilde{\bm{V}} +  2\sum_{j=1}^{n} \boldsymbol{R}_{ij} \boldsymbol{p}_j  \label{7bs}
\end{align}
\end{subequations}
where $\hat{\bm{s}}_{\hat{\imath}} = \sum_{i=1}^{n}\nabla_{\bm{r}_{\hat{\imath}} }(\bm{\lambda}_i^\mathsf{T}\bm{V}_i) $, $\tilde{\bm{V}} = \underline{\bm{V}} -  \boldsymbol{V}_{0}$, $\bm{p}_i = \sum_{\hat{\imath}=1}^{\hat{n}_i} \bm{r}_{\hat{\imath}}$, and $\hat{n}_i$ denotes the number of EVs connected at bus $i$. Note that the exact form of  $\hat{\bm{s}}_{\hat{\imath}}$ is decided based on the bus location of the $\hat{\imath}$th EV, e.g., if the $\hat{\imath}$th EV is connected at bus $k$, then $\hat{\bm{s}}_{\hat{\imath}} = 2\sum_{i=1}^{n} \boldsymbol{R}_{ik}\bm{\lambda}_i$. 

Based on the subgradients in \eqref{7s}, the primal and dual variables can be updated through the PGM by 
\begin{subequations} \label{update}
\begin{align}
\bm{r}_{\hat{\imath}}^{(\ell+1)} &=
\Pi_{\mathcal{R}_{\hat{\imath}}}\left( \bm{r}_{\hat{\imath}}^{(\ell)}-\gamma_{\hat{\imath}} \nabla_{\bm{r}_{\hat{\imath}}} \mathcal{L}\left(\bm{r}^{(\ell)}, \bm{\lambda}^{(\ell)}\right)\right) \label{11a}\\
\bm{\lambda}_i^{(\ell+1)} &=  \Pi_{\mathcal{D}_{i}}\left( 
\bm{\lambda}_i^{(\ell)}+\beta_i \nabla_{\bm{\lambda}_i} \mathcal{L}\left(\bm{r}^{(\ell)}, \bm{\lambda}^{(\ell)}\right)\right)\label{11b}
\end{align}
\end{subequations}
where $\mathcal{D}_{i}= \{\bm{\lambda}_{i} ~|~ \bm{\lambda}_{i} \geq \boldsymbol{0}  \}$ denotes the feasible set of  $\bm{\lambda}_{i}$ and $\beta_i$ denotes the associated dual update step size.

The PGM update in \eqref{update} is scalable \emph{w.r.t.} the number of EVs owing to the parallel computing structure. However, due to the couplings  of decision variables in both the objective function and the global voltage constraint, the primal and dual updates require the exchange of decision variables between all EVs, e.g., calculating the subgradient in \eqref{7as} requires $\bm{r}_{\hat{\imath}}$'s from all EVs. Therefore, without appropriate privacy preservation measures, the inevitable and frequent information exchange can put EVs' private data at breaching risks. To address this concern, we aim to develop a privacy-preserving EV charging control framework via state obfuscation to protect EVs' true decision variables.

\subsection{Privacy-Preserving EV Charging Control Via State Obfuscation}

The goal of privacy preservation is to ensure EV owners' private information is protected during the charging schedules. Specifically, private data of the $\hat{\imath}$th EV are defined to include the charging profiles $\bm{r}_{\hat{\imath}}^{(\ell)}$ in all iterations, charging demand $d_{\hat{\imath}}$, and the maximum charging power $\overline{\bm{r}}_{\hat{\imath}}$. The primal update in \eqref{11a} naturally inherits local privacy preservation owing to the independent projection operation $\Pi_{\mathcal{R}_{i}}$. This is because the private data such as the charging demand $d_{\hat{\imath}}$ and the maximum charging power $\overline{\bm{r}}_{\hat{\imath}}$ are exclusive to the ${\hat{\imath}}$th EV and only contained in the feasible set $\mathcal{R}_{\hat{\imath}}$ for implementing the primal update. Therefore, the local private information is securely retained within $\mathcal{R}_{\hat{\imath}}$ and will not be disclosed to other parties.

Despite the scalability of decentralized  EV charging architectures, they require frequent exchange of EVs' charging profiles through communication channels between EVs and the SO, making the entire system prone to privacy leakages. To resolve this issue, we propose a state-obfuscation-based algorithm that can protect EVs' charging profiles during any planned charging window. 
The basic concept behind state obfuscation is to obfuscate EVs' true decision variables by using the values of random variables drawn from a probability distribution. Regarding a set of mutually independent random variables, e.g., drawn from a normal distribution, we have the following theorem

\noindent \textbf{Theorem 1} \cite{ross1995stochastic}: If $X_1, \ldots, X_z$ are mutually independent normal random variables with means $\mu_1, \ldots, \mu_z$ and variances $\sigma_1^2, \ldots, \sigma_z^2$, then the linear combination $Y = \sum_{i=1}^z c_iX_i$ follows the normal distribution 
$\mathcal{N}(\sum_{i=1}^z c_i \mu_i, \sum_{i=1}^z c_i^2 \sigma_i^2)$. \hfill $\blacksquare$

\begin{figure}[!htb]
    \centering
\includegraphics[width=0.47\textwidth, trim={0.1cm 0cm 0.0cm 0cm},clip]{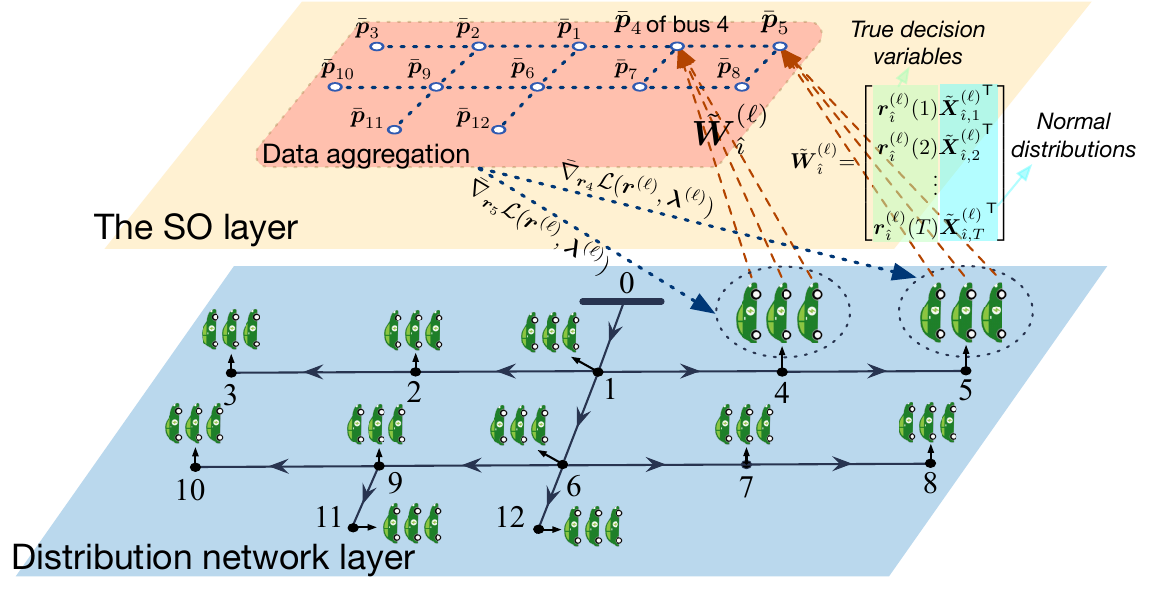}
    \caption{Communicating structure of the proposed privacy-preserving algorithm (the communications of the EVs connected at buses 4 and 5 are given).}
\label{communicating_structure}
     \vspace{-2mm}
\end{figure}

To integrate state obfuscation into EV charging control, we propose a novel communicating architecture, as shown in Fig. \ref{communicating_structure}, for the privacy-preserving algorithm implementation. 
 The EVs in the distribution network layer first send the obfuscated charging profiles to the SO, the SO then aggregates the obfuscated data bus by bus according to EVs' bus locations. Specifically, suppose the $\hat{\imath}$th EV
is connected at bus $i$. At the $\ell$th iteration, the $\hat{\imath}$th EV  uses a random normal variable $X_{\hat{\imath}} \sim \mathcal{N}(\mu_{i}, \sigma_i^2)$
to generate $T$ random sets $\Tilde{X}_{\hat{\imath},t}^{(\ell)} = \{e_{\hat{\imath},t,1}^{(\ell)}, \ldots, e_{\hat{\imath},t,m}^{(\ell)}\}$, $\forall t=1,\ldots,T$ where each set contains $m$ random elements. For clarity, we represent $\Tilde{X}_{\hat{\imath},t}^{(\ell)}$ by a vector as $\Tilde{\bm{X}}_{\hat{\imath},t}^{(\ell)} = [e_{\hat{\imath},t,1}^{(\ell)}, \ldots, e_{\hat{\imath},t,m}^{(\ell)}]^{\mathsf{T}}$.

The $\hat{\imath}$th EV then extracts each element from its decision variable $\bm{r}_{\hat{\imath}}^{(\ell)}$ to calculate
$\bm{r}_{\hat{\imath}}^{(\ell)}(t) \Tilde{\bm{X}}_{\hat{\imath},t}^{(\ell)}$, $\forall t=1,\ldots,T$. Consequently, the $\hat{\imath}$th EV can obtain a total of $T$ new vectors and reformulate them into ${\Tilde{\bm{W}}_{\hat{\imath}}}^{(\ell)}$, defined by 
\begin{equation} \label{11st}
\Tilde{\bm{W}}_{\hat{\imath}}^{(\ell)} =  [\bm{r}_{\hat{\imath}}^{(\ell)}(1) {{}\Tilde{\bm{X}}_{\hat{\imath},1}^{(\ell)}}^{\mathsf{T}}, \ldots, \bm{r}_{\hat{\imath}}^{(\ell)}(T) {{}\Tilde{\bm{X}}_{\hat{\imath},T}^{(\ell)}}^{\mathsf{T}} ]^{\mathsf{T}}.
\end{equation}

Subsequently, instead of sending the true decision variables directly, all EVs send their $\tilde{\bm{W}}_{\hat{\imath}}^{(\ell)}$'s, i.e., the obfuscated states, to the SO. Then the SO computes the sum of the received obfuscated states using
\begin{equation} \label{12sss}
    \bm{Y}_i^{(\ell)} = \sum_{\hat{\imath}=1}^{\hat{n}_i}   \Tilde{\bm{W}}_{\hat{\imath}}^{(\ell)}
\end{equation}
for the EVs connected at the $i$th bus.

As shown in \eqref{11st}, every element in $\bm{r}_{\hat{\imath}}^{(\ell)}$ is obfuscated and expanded  by $m$ random values. To retrieve the summed charging profiles for all EVs connected at bus $i$, the SO needs to calculate the mean of  every $m$ elements in $\bm{Y}_i^{(\ell)}$. For example, for the first $m$ elements, the SO calculates
$(\sum_{\kappa=1}^{m} \bm{Y}_i^{(\ell)}(\kappa))/m$ that is equal to $\sum_{\hat{\imath}=1}^{\hat{n}_i} \bm{r}_{\hat{\imath}}^{(\ell)} (1) \bar{\mu}_{i}$ where $\bar{\mu}_{i}$ is an approximation of $\mu_{i}$. Suppose the SO knows the true mean $\mu_{i}$, the SO can acquire  $\sum_{\hat{\imath}=1}^{\hat{n}_i}  \bm{r}_{\hat{\imath}}^{(\ell)} \bar{\mu}_{i}$, and further obtain $\tau\sum_{\hat{\imath}=1}^{\hat{n}_i}  \bm{r}_{\hat{\imath}}^{(\ell)}$ where $\tau  = \bar{\mu}_{i}/\mu_{i}$ denotes the approximation error. Therefore, the SO now has the approximated active power consumption $\bar{\bm{p}}_i^{(\ell)} = \tau \sum_{\hat{\imath}=1}^{\hat{n}_i}  \bm{r}_{\hat{\imath}}^{(\ell)}$ of bus $i$, and it repeats the procedure to obtain $\bar{\bm{p}}_i^{(\ell)}$, $\forall i=1,\ldots,n$. 

Finally, the SO estimates the subgradients in \eqref{7s} using the approximated active power consumption, then broadcasts \eqref{7as} to the $i$th EV while utilizing \eqref{7bs} to conduct dual updates. Thereafter, EV $i$ can update its decision variable $\bm{r}_{\hat{\imath}}^{(\ell)}$ in parallel using \eqref{11a}.

The step-by-step process of the proposed approach is outlined in 
\textbf{Algorithm \ref{alg_1}}.

\begin{algorithm}
\caption{Decentralized privacy-preserving EV charging control  via state obfuscation}
\begin{algorithmic}[1]

\State EVs initialize decision variables, tolerance $\epsilon_0$,  iteration counter $\ell=0$, and maximum iteration $\ell_{max}$.

\While{ $\epsilon_{\hat{\imath}}^{(\ell)} > \epsilon_0$ and $\ell < \ell_{max}$}

\State The $\hat{\imath}$th EV connected at bus $i$ generates  a normal random variable $X_{\hat{\imath}} \sim \mathcal{N}(\mu_{i}, \sigma_i^2)$ and draws random elements from $X_{\hat{\imath}}$  to obtain $\Tilde{\bm{X}}_{\hat{\imath},t}^{(\ell)} = [e_{\hat{\imath},t,1}^{(\ell)}, \ldots, e_{\hat{\imath},t,m}^{(\ell)}]^{\mathsf{T}}$, $\forall t=1,\ldots,T$. 

\State The $\hat{\imath}$th EV uses the elements of $\bm{r}_{\hat{\imath}}^{(\ell)}$ to calculate
$\bm{r}_{\hat{\imath}}^{(\ell)}(t) \Tilde{\bm{X}}_{\hat{\imath},t}^{(\ell)}$, $\forall t=1,\ldots,T$ elementwisely. Then each EV formulates $\Tilde{\bm{W}}_{\hat{\imath}}^{(\ell)}$ and sends it to the SO. 

\State The SO calculates the summation $\bm{Y}_i^{(\ell)}$ using \eqref{12sss} for each bus, then calculates the mean of every $m$ elements in $\bm{Y}_i^{(\ell)}$, to obtain the approximated $\bar{\bm{p}}_i^{(\ell)}$, $\forall i =1, \ldots, n$.

\State The SO estimates the subgradient in  \eqref{7as} and broadcasts it to the $\hat{\imath}$th EV.

\State  The $\hat{\imath}$th EV updates $\bm{r}_{\hat{\imath}}^{(\ell)} \rightarrow \bm{r}_{\hat{\imath}}^{(\ell+1)}$ by PGM using \eqref{11a}, then calculates the error $\epsilon_{\hat{\imath}}^{(\ell)}$.

\State  The SO updates the dual variables $\bm{\lambda}_{i}^{(\ell)} \rightarrow \bm{\lambda}_{i}^{(\ell+1)}$, $\forall i=1,\ldots,n$ using \eqref{11b}. 

\State  $\ell=\ell+1$. 

\EndWhile
\end{algorithmic}
\label{alg_1}
\end{algorithm}

\noindent \textbf{Theorem 2}: \textbf{Algorithm \ref{alg_1}} has an accuracy level of $\tau$. With appropriate choices of $\sigma$ and $m$, the convergence of primal and dual variables is guaranteed. 
\hfill$\blacksquare$

\textbf{Theorem 2} states the correctness and convergence of the proposed algorithm. By carrying out \textbf{Algorithm 1}, the subgradients in \eqref{7s} can be efficiently approximated and calculated since the mean of $\bm{Y}_i^{(\ell)}$ can be used to retrieve $\bar{\bm{p}}_i^{(\ell)}$ that is an estimation of $\bm{p}_i^{(\ell)}$. When determining the accuracy level, the standard error of the mean (SEM), defined by $SEM_{m} = \sigma/\sqrt{m}$, can quantify how a larger sample size produces more precise estimates of the means.


\noindent \textbf{Remark 1 :} Without the loss of generality,  $\mu_{i}$ was set uniformly across all EVs connected at the same bus to avoid over-complicated algorithm implementation. In a broader scenario, the mean values of different EVs can be chosen independently. In other words, the mean value $\mu_{\hat{\imath}}$ serves as a unique key between the $\hat{\imath}$th EV and the SO. \hfill $\square$

\section{Privacy Analysis}

\subsection{Privacy and Attack Models}


To preserve EV owner's privacy, two types of adversaries are considered: 1) An \emph{honest-but-curious adversary} is an agent who adheres to the algorithm but intends to utilize the accessible data to infer private information of other participants, and 2) an \emph{external eavesdropper} is an external attacker who wiretaps communication links to obtain the private information of the participants.

\subsection{Privacy Analysis}

\textbf{Algorithm \ref{alg_1}} allows EVs to use the values of random variables drawn from a normal distribution to protect the true decision variables $\bm{r}_{\hat{\imath}}^{(\ell)}$'s. The privacy preservation properties of \textbf{Algorithm \ref{alg_1}} are given by the following theorem

\noindent \textbf{Theorem 3}: \textbf{Algorithm \ref{alg_1}} preserves the private data of EV owners against both honest-but-curious adversaries and external eavesdroppers. \hfill$\blacksquare$

\textit{Proof}: Proof of \textbf{Theorem 3} is approached from the adversaries' perspective based on the data they can access. From the view of an honest-but-curious adversary, suppose both EV $\hat{\imath}_1$ and EV $\hat{\imath}_2$ are connected at the same bus, and EV $\hat{\imath}_1$ is curious in inferring the charging profiles of EV $\hat{\imath}_2$. At the $\ell$th iteration, EV $\hat{\imath}_1$ can have access to the data set $\mathcal{A}_{\hat{\imath}_1}^{(\ell)} = \{ \bm{r}_{\hat{\imath}_1}^{(\ell)}, \overline{\bm{r}}_{\hat{\imath}_1}, d_{\hat{\imath}_1},\mathcal{R}_{\hat{\imath}_1}, \gamma_{\hat{\imath}}, \bar{\nabla}_{\bm{r}_{\hat{\imath}_1} } \mathcal{L}(\bm{r}^{(\ell)},\bm{\lambda}^{(\ell)} ) \}$ where $\overline{\bm{r}}_{\hat{\imath}_1}$,  $d_{\hat{\imath}_1}$, $\mathcal{R}_{\hat{\imath}_1}$, and $\gamma_{\hat{\imath}}$  are private information of EV $\hat{\imath}_{\hat{\imath}_1}$ and kept to EV $\hat{\imath}_1$ locally. The local information, therefore, cannot provide any useful information in inferring $\bm{r}_{\hat{\imath}_2}^{(\ell)}$. Besides, EV $\hat{\imath}_1$ also has access to the approximated subgradient $\bar{\nabla}_{\bm{r}_{\hat{\imath}_1} } \mathcal{L}(\bm{r}^{(\ell)},\bm{\lambda}^{(\ell)})$ that is calculated by the SO. However, the baseline load $\bm{p}_b$ and the adjacency matrix $\bm{R}$ are held by the SO, and therefore remain invisible to any EVs. Therefore, EV $\hat{\imath}_1$ cannot infer the charging profiles $\bm{r}_{\hat{\imath}_2}$ of EV $\hat{\imath}_2$ based on its accessible information contained in  $\mathcal{A}_{\hat{\imath}_1}^{(\ell)}$. 

For any external eavesdropper, by wiretapping the communication channels at the $\ell$th iteration, it can obtain the information set $\mathcal{E} = \{ \Tilde{\bm{W}}_{\hat{\imath}}^{(\ell)}, \bar{\nabla}_{\bm{r}_{\hat{\imath}} } \mathcal{L}(\bm{r}^{(\ell)},\bm{\lambda}^{(\ell)}), \forall \hat{\imath} =1, \ldots, \hat{n} \}$. Suppose an external eavesdropper knows the protocols of \textbf{Algorithm \ref{alg_1}}. To infer $\bm{r}_{\hat{\imath}}^{(\ell)}$ by using $\mathcal{E}$, it still needs to know the cardinality $m$ and the mean value $\mu_{i}$ that is associated with the random variable $X_{\hat{\imath}}$. Though the approximated subgradient $\bar{\nabla}_{\bm{r}_{\hat{\imath}} } \mathcal{L}(\bm{r}^{(\ell)},\bm{\lambda}^{(\ell)})$ could potentially reveal the converging direction of the decision variable, the 
eavesdropper is still blind from $\gamma_{\hat{\imath}}$
and $\mathcal{R}_{\hat{\imath}}$, therefore unable to imitate the primal update in \eqref{11a}. \hfill$\square$

\noindent \textbf{Remark 2:} A trade-off between the level of security and computing cost exists in  \textbf{Algorithm \ref{alg_1}}. When a specific accuracy requirement is decided by $\tau$ under a fixed sample size $m$, a smaller variance $\sigma_{i}^2$ will result in a smaller SEM and therefore require fewer data points to achieve the accuracy standard. The proposed state obfuscation refines $k$-anonymity \cite{samarati1998protecting} by introducing randomization of $m$ anonymous random variables for each true value. Though a smaller variance would result in less computation and communication cost, it can also lead to a higher degree of similarity in $\Tilde{\bm{X}}_{\hat{\imath},t}^{(\ell)}$, thus compromising the level of randomization and privacy.
\hfill $\square$

\section{Simulation Results}

The effectiveness of the proposed obfuscation-based privacy-preserving EV control strategy is verified through the simplified single-phase IEEE 13-bus test feeder as shown in Fig. \ref{13_distribution_network}.  The baseline load profile was taken and scaled from California Independent System Operator on 09/16/2021 and 09/17/2021 \cite{CISO}. We consider the penetration level of 7 EVs per bus, and in total 84 EVs are connected to the distribution network. The charging demands of all EVs randomly distribute in $[10,40]$ kWh. The maximum charging power $\bar{r}_{\hat{\imath}}$'s are uniformly set to be 6.6 kW based on the level-2 EV charging standards and the charging efficiency is set to be $\eta = 0.85$. The valley-filling horizon is set to begin at 19:00 and lasts until 7:00 the next morning. The entire control horizon is divided into $T=48$ time slots with 15-minute resolution. It is required that, by the end of the valley-filling period, all EVs need to be charged to the desired energy levels. The primal update step sizes are chosen based on experience as $\gamma_{\hat{\imath}} = 4\times10^{-4}$, $\forall \hat{\imath} =1,\dots,\hat{n}$, and the dual update step sizes are $\beta_i = 2\times10^{-3}$, $\forall i =1,\dots, n$. Initial values of $\bm{r}_{\hat{\imath}}^{(0)}$'s and $\bm{\lambda}_i^{(0)}$'s are all set to be zeros. The normal random variables $\Tilde{X}_{\hat{\imath}}$'s generated by EVs connected at bus $i$ follow the nomral distribution $\Tilde{X}_{\hat{\imath}} \sim \mathcal{N}(\mu_{i} = 1, \sigma_{i}^2 = 0.2)$. The cardinality of $\Tilde{\bm{X}}_{\hat{\imath}}$'s is set uniformly to be $m=40$. 

\begin{figure}[h!]
\vspace{-2mm}
\centering
\begin{minipage}[t]{.24\textwidth}
  \centering
\includegraphics[width=1\textwidth, trim={0cm 0cm 0cm 0.23cm},clip]{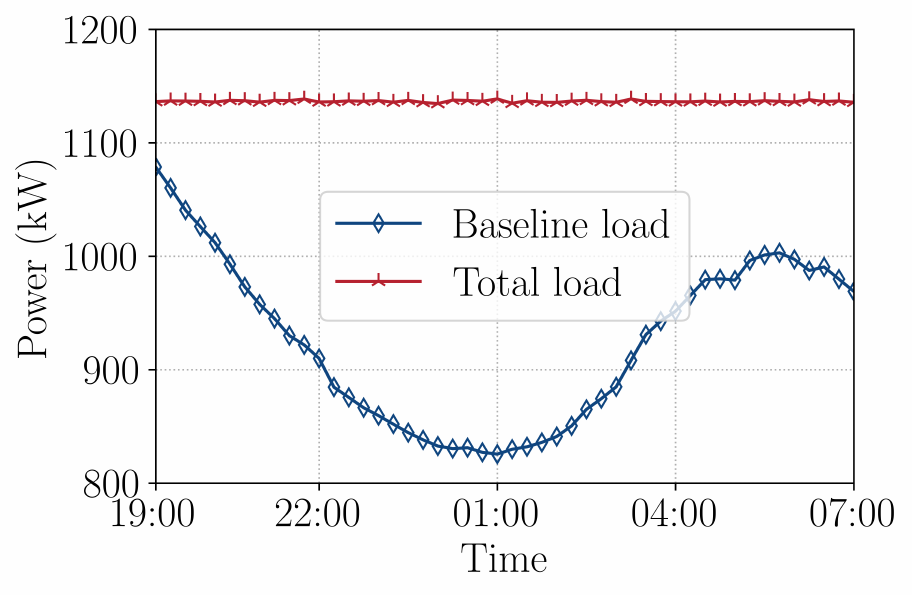}
  \caption{The overnight valley filling.}
  \label{valley_filling}
\end{minipage}%
\hfill
\begin{minipage}[t]{.225\textwidth}
  \centering
  \includegraphics[width=1\textwidth, trim={0cm 0cm 0cm 0cm},clip]
{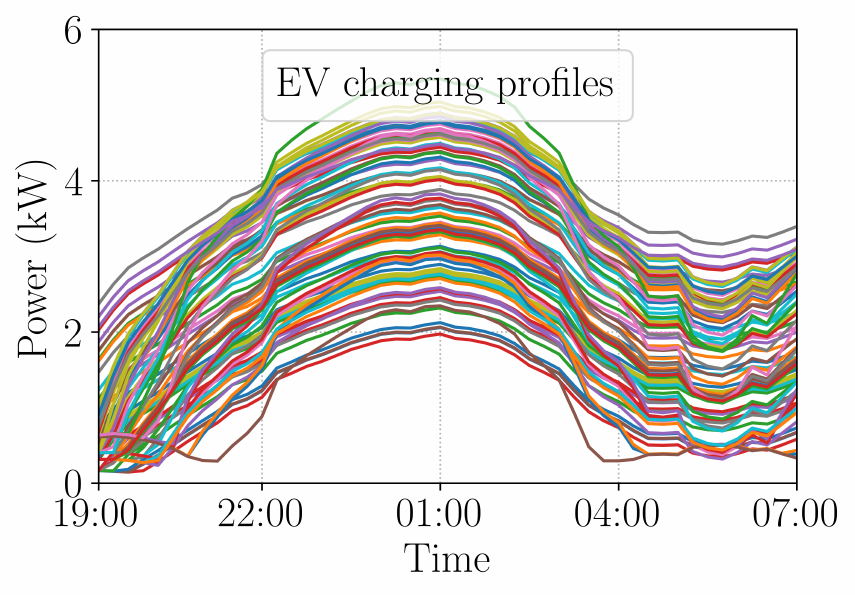}
  \caption{Optimal charging profiles of all EVs.}
  \label{charging_profiles}
\end{minipage}
\vspace{-3mm}
\end{figure}

By applying \textbf{Algorithm \ref{alg_1}}, 
Fig. \ref{valley_filling}  shows that the baseline load was flattened by using EVs' charging load. The optimal charging profiles of all EVs are shown in Fig. \ref{charging_profiles}. At around 1:00 a.m., when the baseline load reaches its minimum, all EVs charge at their highest power.



\begin{figure}[!htb]
 \vspace{-2mm}
  \centering
\includegraphics[width=0.44\textwidth,trim={0cm 0cm 0cm 0cm},clip]{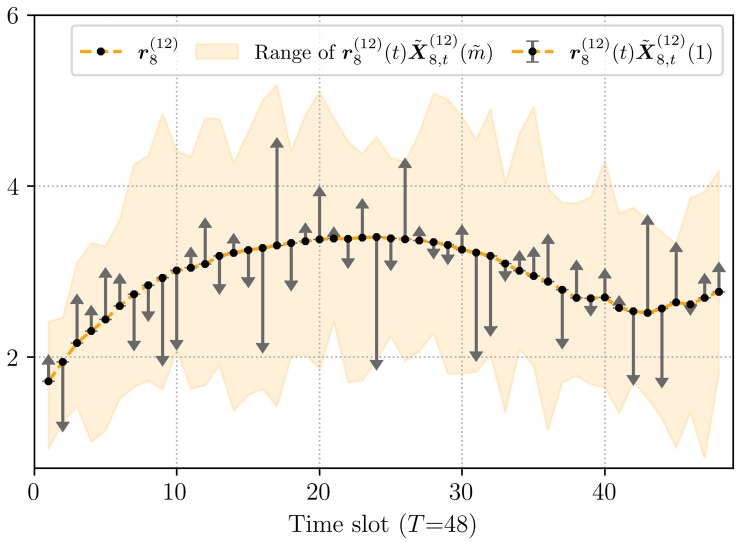}%
  \caption{The true and obfuscated data generated by EV 8 at the 12th iteration.}
  \label{fig6}
  \vspace{-5mm}
\end{figure}

\begin{figure}[!htb]
\centering\includegraphics[width=0.48\textwidth, trim={0cm 0cm 0cm 0cm},clip]{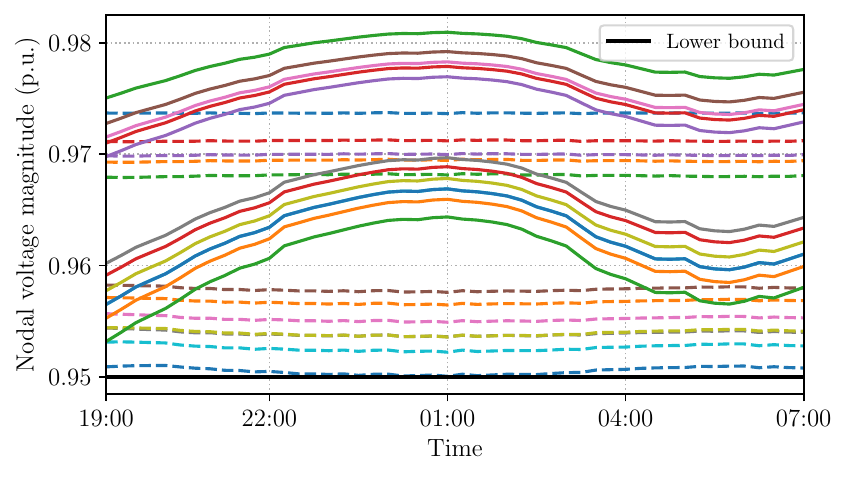}
    \caption{Nodal voltage magnitudes of 12 buses under baseline load (solid lines) and total load (dashed lines).}
\label{voltage}
     \vspace{-2mm}
\end{figure}

To observe the privacy features, Fig. \ref{fig6} presents the random values $\bm{r}_{8}^{(12)}(t) \Tilde{\bm{X}}_{8,t}^{(12)}$, $\forall t=1,\ldots,T$ that were generated by EV 8 at the 12th iteration. The true charging profile $\bm{r}_{8}^{(12)}$ was obfuscated into $\bm{r}_{8}^{(12)}(t) \Tilde{\bm{X}}_{8,t}^{(12)}(\Tilde{m})$,  $\forall t=1,\ldots, T, \Tilde{m} = 1,\ldots,m$. The range of the obfuscated data is shown by the shaded area, where the obfuscation achieves nearly 50\% randomization of the original data. 
Fig. \ref{voltage}  shows the nodal voltage magnitudes of 12 buses in the distribution network, where all voltage magnitudes are above the lower voltage limit (0.95 p.u., the black line).


\section{Conclusion}

In this paper, we proposed a novel privacy-preserving decentralized algorithm to achieve privacy preservation and scalability in large-scale multi-agent cooperative optimization, particularly in the context of cooperative EV charging control. The proposed algorithm enables EVs to protect their decision variables via state obfuscation while facilitating the cooperation between EVs and the SO to achieve overnight valley filling. The privacy guarantees were theoretically analyzed and evaluated against honest-but-curious adversaries and external eavesdroppers. Simulations on an EV charging control problem validated the accuracy, efficiency, and privacy preservation properties of the proposed approach. 

 
\bibliographystyle{IEEEtran}

\bibliography{bibliography}

\end{document}